%% file: yama.tex
\documentclass[11pt]{article}
\input amssym.def
\input amssym.tex
\usepackage{epsfig,psfig}
\newtheorem{theorem}{Theorem}

\def\binom#1#2{{#1}\choose{#2}}

\def\slfrac#1#2{\hbox{\kern.1em %
 \raise.5ex\hbox{\the\scriptfont0 #1}\kern-.11em %
 /\kern-.15em\lower.25ex\hbox{\the\scriptfont0 #2}}}

\newcommand{\eqn}[1]{(\ref{#1})}
\newcommand{\hsp}{\hspace*{\parindent}}

\newcommand{\eeq}{\end{equation}}
\newcommand{\beql}[1]{\begin{equation}\label{#1}}

\newcommand{\sO}{{\cal O}}
\newcommand{\sF}{{\cal F}}

\newcommand{\sI}{{\cal I}}

\newcommand{\sA}{{\cal A}}

\newcommand{\ZZ}{{\Bbb Z}}
\newcommand{\RR}{{\Bbb R}}
\newcommand{\FF}{{\Bbb F}}

\newcommand{\la}{\lambda}

\newcommand{\Om}{\Omega}

\newcommand{\dd}{\ldots}

\newcommand{\df}{\displaystyle\frac}

\newcommand{\sG}{{\cal G}}

\newcommand{\sD}{{\cal D}}

\newcommand{\sC}{{\cal C}}

\newcommand{\sP}{{\cal P}}
\newcommand{\sQ}{{\cal Q}}
\newcommand{\sT}{{\cal T}}

\makeatletter
\def\@sect#1#2#3#4#5#6[#7]#8{\ifnum #2>\c@secnumdepth
     \def\@svsec{}\else
     \refstepcounter{#1}\edef\@svsec{\csname the#1\endcsname.\hskip .75em }\fi
     \@tempskipa #5\relax
      \ifdim \@tempskipa>\z@
        \begingroup #6\relax
          \@hangfrom{\hskip #3\relax\@svsec}{\interlinepenalty \@M #8\par}%
        \endgroup
       \csname #1mark\endcsname{#7}\addcontentsline
         {toc}{#1}{\ifnum #2>\c@secnumdepth \else
                      \protect\numberline{\csname the#1\endcsname}\fi
                    #7}\else
        \def\@svsechd{#6\hskip #3\@svsec #8\csname #1mark\endcsname
                      {#7}\addcontentsline
                           {toc}{#1}{\ifnum #2>\c@secnumdepth \else
                             \protect\numberline{\csname the#1\endcsname}\fi
                       #7}}\fi
     \@xsect{#5}}
\def\@begintheorem#1#2{\it \trivlist \item[\hskip \labelsep{\bf #1\ #2.}]}

\def\plain{plain}\ifx\fmtname\plain\csname fi\endcsname
     
     \input docstrip
     \preamble

     Do not distribute the stripped version of this file.
     The checksum in the header refers to the documented version.

     \endpreamble
     \generateFile{here.sty}{t}{\from{here.doc}{}}
     \endinput
\fi
\ifcat a\noexpand @\let\next\relax\else\def\next{%
    \documentstyle[here,doc]{article}\MakePercentIgnore}\fi\next
\ifx\@Hxfloat\@Hundef\else\expandafter\endinput\fi
\let\@Hxfloat\@xfloat
\def\@xfloat#1[{\@ifnextchar{H}{\@HHfloat{#1}[}{\@Hxfloat{#1}[}}
\def\@HHfloat#1[H]{%
\expandafter\let\csname end#1\endcsname\end@Hfloat
\vskip\intextsep\vbox\bgroup\def\@captype{#1}\parindent\z@
\ignorespaces}
\def\end@Hfloat{\egroup\vskip \intextsep}

\makeatother

\begin{document}
\begin{center}
{\Large {\bf Packing Planes in Four Dimensions and Other Mysteries}} \\
\vspace{1.5\baselineskip}
{\em N. J. A. Sloane} \\
\vspace*{.5\baselineskip}
Information Sciences Research \\
 AT\&T Labs-Research \\
180 Park Avenue, Room C233 \\
Florham Park, NJ 07932-0971 \\
\vspace{1.5\baselineskip}
April 7, 1998 \\
\vspace{1.5\baselineskip}
{\bf ABSTRACT}
\vspace{.5\baselineskip}
\end{center}
\setlength{\baselineskip}{1.5\baselineskip}

How should you choose a good set of (say) 48 planes in four dimensions?
More generally, how do you find packings in Grassmannian spaces?
In this article I give a brief introduction to the work that I have been doing on this problem in collaboration with A. R. Calderbank, J. H. Conway,
R. H. Hardin, E. M. Rains and P. W. Shor.
We have found many nice examples of specific packings (70 4-spaces in 8-space,
for instance), several general constructions, and an
embedding theorem which shows that a packing in Grassmannian space $G(m,n)$
is a subset of a sphere in $\RR^D$,
$D= (m+2) (m-1)/2$, and leads to a proof that many of our packings are optimal.
There are a number of interesting unsolved problems.

\section{Introduction}
\hsp
In my talk at the Yamagata conference on ``Algebraic Combinatorics and Related Topics'' (November 1997) I discussed two problems, (a)~finding packings in Grassmannian manifolds, and (b)~constructing error-correcting codes for quantum computation, and tried to show how our purely numerical investigations into the first problem had led to theoretical advances in both subjects.
This work has been presented in a series of papers,
\cite{grass3},
\cite{CRSS96}, \cite{qc2},
\cite{grass}, \cite{grass2}.
In the present paper I will give a brief introduction to our work on the first problem, referring the reader to the above references for further information.

For a long time I have been interested in various kinds of packing problems:
packings in Hamming space (i.e. error-correcting codes \cite{MS77}), in
Euclidean space (i.e. the sphere-packing problem \cite{SPLAG}), or on the
sphere (i.e. spherical codes \cite[Chap. 3]{SPLAG}, \cite{HSS94},
\cite{HSS96}).
R. H. Hardin and I have had some success in using numerical optimization techniques to search for spherical codes and other kinds of geometrical designs.
This work is described in the papers
\cite{HS92}, \cite{HS93}, \cite{Me174}, \cite{Me189}, \cite{Me203}, \cite{Me204} (others are in preparation).
Figure 1 (on the next page) shows an example
that arose from our search for good coverings, i.e. sets of points with small covering radius.
The figure shows 48002 points on the sphere with average angular separation of about 1 degree
(the minimal angular separation is .802 degrees,
the maximal angle is 1.077, the average angle is .99948 and the standard deviation is .047).

A few years ago, a statistician, Dianne Cook, asked if we could apply the same techniques to find packings of (say) 48 planes in 4-dimensional Euclidean space $\RR^4$.
More generally, for given values of $N$, $n$ and $m$,
how should one arrange $N$ $n$-dimensional subspaces of
$\RR^m$ (all passing through the origin) so that they are as far apart as
possible?

The rest of the paper will describe what happened when we tried to attack this problem.

\section{Packings in Grassmannian Space}
\hsp
The {\em Grassmannian space} $G(m,n)$ is the set of all $n$-dimensional
subspaces of real Euclidean $m$-dimensional space $\RR^m$.
This is a homogeneous space isomorphic to $O(m) / (O(n) \times O(m-n))$, and forms a compact Riemannian manifold of dimension $n(m-n)$.

Before we can investigate packings in this space, we must decide how to measure
the distance between two $n$-spaces
$P$, $Q \in G(m,n)$.
The {\em principal angles} $\theta_1, \ldots , \theta_n \in [0, \pi/2 ]$ between
$P$ and $Q$ are defined by (we follow \cite{GVL89}, p.~584)
$$\cos \theta_i = \max_{u \in P} \max_{v \in Q} u \cdot v = u_i \cdot v_i ~,$$

\clearpage
\begin{figure}[htb]
\centerline{\psfig{file=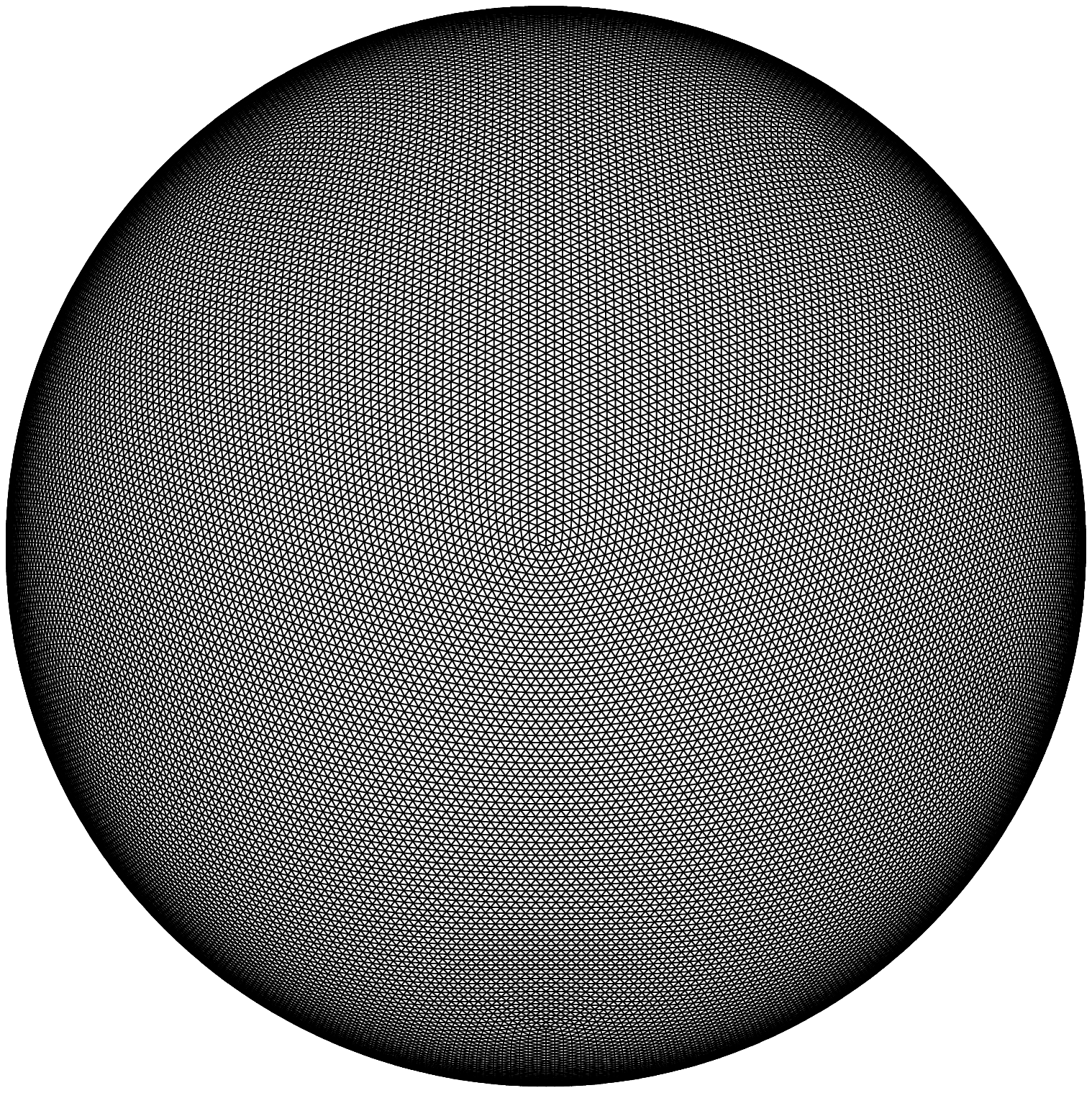}}
\end{figure}
\clearpage
\noindent
for $i=1, \ldots , n$,
subject to $u \cdot u = v \cdot v =1$, $u \cdot u_j =0$, $v \cdot v_j =0$ $(1 \le j \le i-1)$.
The vectors $\{u_i\}$ and $\{v_j \}$ are {\em principal vectors}
corresponding to the pair $P$ and $Q$.

Wong \cite{Won67} shows
that the {\em geodesic distance} on $G(m,n)$ between $P$ and $Q$ is\footnote{The geodesic distance is unique except for the single case of $G(4,2)$.}
\beql{Eq1}
d_g (P,Q) = \sqrt{\theta_1^2 + \cdots + \theta_n^2} ~.
\eeq
However, this definition has one drawback:
it is not everywhere differentiable.
Consider the case $n=1$, for example, and hold one line $P$ fixed while rotating another line $Q$ (both passing through the origin).
As the angle $\phi$ between $P$ and $Q$ increases from 0 to $\pi$,
the principal angle $\theta_1$
increases from 0 to $\pi /2$ and then falls to 0, and is non-differentiable at $\pi /2$ (see Fig.~\ref{fg1}).
\setcounter{figure}{1}
\begin{figure}[htb]
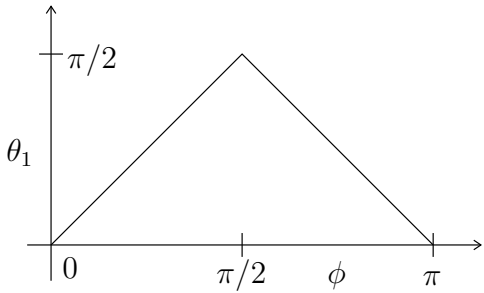

\begin{center}
\input tent_text.tex
\end{center}
\caption{Principal angle $\theta_1$ between two lines as the angle between them increases from 0 to $\pi$.}
\label{fg1}
\end{figure}

Although one might expect this non-differentiability to be a mere technicality, it does
in fact cause considerable difficulties for our optimizer, especially
in higher dimensions in cases when many distances fall in the neighborhood of singular
points of $d_g$.

An alternative measure of distance, which we call the {\em chordal distance}, is given by
\beql{Eq2}
d_c (P,Q) = \sqrt{\sin^2 \theta_1 + \cdots + \sin^2 \theta_n} ~.
\eeq
The reason for the name will be revealed later.
This approximates the geodesic distance when the subspaces are close,
has the property that its square is differentiable everywhere,
and, as we shall attempt to demonstrate,
has a number of other desirable features.

A third definition has been used by Asimov \cite{Asi85} and Golub and Van~Loan \cite{GVL89}, p.~584,
namely
$$d_m (P,Q) = \max_{i=1, \dd, n} \theta_i ~.$$
This shares the vices of the geodesic distance.

Of course for $n=1$ all three definitions are equivalent, in the sense that
they lead to the same optimal packings.

We can now state the packing problem:
given $N,n,m$, find a set of $n$-spaces $P_1, \dd, P_N \in G(m,n)$ so that
$\min\limits_{i \neq j} d(P_i, P_j)$ is as large as possible,
where $d$ is either geodesic or chordal distance.
Since $G(m,n)$ is compact, the problem is well-defined.
Because $G(m,n)$ and $G(m, m-n)$ are essentially the same space, we may assume $n \le m/2$.

We also need some further terminology.
A {\em generator matrix} for an $n$-space $P\in G (m,n)$ is an $n \times m$ matrix whose rows span $P$.
The orthogonal group $O(m)$ acts on $G(m,n)$ by right multiplication of generator matrices.
The {\em automorphism group} of a subset $\{P_1, \dd, P_N \} \subset G(m,n)$ is the subset of $O(m)$ which fixes or permutes these planes.

By applying a suitable element of $O(m)$ and choosing appropriate
basis vectors for the spaces, any given pair of $n$-spaces $P,Q$ with $n \le m/2$ can be assumed to have generator matrices
\beql{Eq3}
\left[
\matrix{
1 & 0 & \cdots & 0 & 0 & 0 & \cdots & 0 & 0 & \cdots & 0 \cr
0 & 1 & \cdots & 0 & 0 & 0 & \cdots & 0 & 0 & \cdots & 0 \cr
\cdot & \cdot & \cdots & \cdot & \cdot & \cdot & \cdots & \cdot & \cdot & \cdots & \cdot \cr
0 & 0 & \cdots & 1 & 0 & 0 & \cdots & 0 & 0 & \cdots & 0\cr}
\right]
\eeq
and
\beql{Eq4}
\left[
\matrix{\cos \theta_1 & 0 & \cdots & 0 & \sin \theta_1 & 0 & \cdots & 0 & 0 & \cdots & 0 \cr
0 & \cos \theta_2 & \cdots & 0 & 0 & \sin \theta_2 & \cdots & 0 & 0 & \cdots & 0 \cr
\cdot & \cdot & \cdots & \cdot & \cdot & \cdot & \cdots & \cdot & \cdot & \cdots & \cdot \cr0 & 0 & \cdots & \cos \theta_n & 0 & 0 & \cdots & \sin \theta_n & 0 & \cdots & 0 \cr}
\right]
\eeq
respectively, where $\theta_1 , \dd, \theta_n$ are the principal angles (\cite{Won67}, Theorem~2).

\section{Packing lines in $\RR^m$}
\hsp
The next step was to apply a suitably modified version of our optimizer in order to find a collection of good packings.
The basic algorithm is described in \cite{HS93} and we shall not say much about it here.
We had to overcome various obstacles caused by the nondifferentiability of the geodesic distance.
Whereas we could compute partial derivatives of the chordal distance analytically, for the geodesic distance we used numerical differentiation.

Our initial computations concerned packings in $G(m,1)$;
that is, packings of lines through the origin in $\RR^m$, or equivalently the problem of packing points on the sphere in $\RR^m$ so that if $P$ is a point in
the packing, so is $-P$.
Such sets of points form {\em antipodal spherical codes}.
This is a classical problem, of course, that is the subject of many papers.
The 3-dimensional problem is discussed for example in
\cite{FT65}, \cite{Ros93}, \cite{Ros94}, \cite{Ros96},
\cite{Rut45}, \cite{SchW51}, \cite{Stroh63}, \cite{TaG83}.

Nevertheless we were able to find many packings that were better than those previously known, even in three dimensions.
We refer the reader to \cite{grass}, \cite{HSS94} for the details.
Here we will just give one table, showing
the minimal angle $\theta_1$ of the best packing we have found of $N$ lines,
and for comparison the minimal angle $\theta'_1$ of the best packing known of $2N$ points
on $S^2$ (taken from \cite{HSS94}).
We see that requiring a packing of $2N$ points on $S^2$ to be antipodal is a definite handicap:
only in the cases of 6 and 12 points do the antipodal and unrestricted packings coincide.
Decimals in the tables have been rounded to four decimal places.

The last two columns of the table specify the largest automorphism group\footnote{That is, the subgroup of $O(3)$ that fixes or permutes the $2N$ points.} we have found of any such best antipodal packing of $2N$ points.
The fourth column gives the order of the group and its name both in the orbifold notation
(cf. \cite{CS96}) and as the double cover of a rotation group.
The symbol $\pm \sG$ indicates that the group consists of the matrices
$\pm M$ for $M \in \sG$, where $\sG$ is a cyclic $( \sC )$, dihedral $( \sD )$,
tetrahedral $(\sT)$,
octahedral $(\sO)$ or icosahedral $(\sI)$ group.
In each case the subscript gives the order of the rotation group.
In some cases the best packings can be obtained by taking the diameters of a known
polyhedron, and if so this is indicated in the final column of the table.

The entries for $N \le 6$ were shown to be optimal by Fejes T\'{o}th in 1965 \cite{FT65}
(see also Rosenfeld \cite{Ros94}),
and the 7-line arrangement will be proved optimal in Section~5.
The solutions for $N \ge 8$ are the best found with over 15000
random starts with our optimizer.
There is no guarantee that these are optimal, but experience with similar problems
suggests that they will be hard to beat and in any case will be not far from optimal.

For $N=1,2,3,6$ the solutions are known to be unique,
for $N=4$ there are precisely two solutions (\cite{FT65}, \cite{Ros93}, \cite{Ros94}), and for $N=5,7,8$ the solutions appear to be unique.
For larger values of $N$, however, the solutions are often not unique.
For $N=9$ there are two different solutions, and in the range $N \le 30$ the solutions
for 10, 22, 25, 27, 29 lines (and possibly others)
contain lines that ``rattle'', that is, lines which can be moved freely over a small range of angles
without affecting the minimal angle.
The table only goes as far as $N=28$ lines.
For more a extensive table in $G(3,1)$ and for tables of packings in $G(m,1)$ for $m > 3$, see \cite{grass}.

\begin{table}[htb]
\caption{Best packings known of $N$ lines in $\RR^3$ ($\theta_1$ gives angle).}

$$
\begin{array}{|rccr@{}ll|} \hline
N & \min ~\theta_1 & \min ~\theta'_1 & \multicolumn{2}{c}{\mbox{group}} & \mbox{polyhedron} \\ \hline
2~~~ & 90.0000 & 109.4712 & *224 = & \pm \sD_8 & \mbox{square} \\
3~~~ & 90.0000 & 90.0000 & *432 = & \pm \sO_{24} & \mbox{octahedron} \\
4~~~ & 70.5288 & 74.8585 & *432 = & \pm \sO_{24} & \mbox{cube} \\
5~~~ & 63.4349 & 66.1468 & 2\! * \!5 = & \pm \sD_{10} & \mbox{pentagonal antiprism} \\
6~~~ & 63.4349 & 63.4349 & *532 = & \pm \sI_{60} & \mbox{icosahedron} \\
7~~~ & 54.7356 & 55.6706 & *432 = & \pm \sO_{24} & \mbox{rhombic dodecahedron} \\
8~~~ & 49.6399 & 52.2444 & \times = & \pm \sC_1 & ~ \\
9~~~ & 47.9821 & 49.5567 & 3\times = & \pm \sC_3 & ~ \\
9~~~ & 47.9821 & 49.5567 & 2* = & \pm \sC_2 & ~ \\ 
10~~~ & 46.6746 & 47.4310 & *226 = & \pm \sD_{12} & \mbox{hexakis bi-antiprism}\\
11~~~ & 44.4031 & 44.7402 & 2\! * \! 5 = & \pm \sD_{10} & ~ \\
12~~~ & 41.8820 & 43.6908 & *432 = & \pm \sO_{24} & \mbox{rhombicuboctahedron} \\
13~~~ & 39.8131 & 41.0377 & 2* = & \pm \sC_2 & ~ \\
14~~~ & 38.6824 & 39.3551 & \times = & \pm \sC_1 & ~  \\
15~~~ & 38.1349 & 38.5971 & 2 \!* \!5 = & \pm \sD_{10} & ~ \\
16~~~ & 37.3774 & 37.4752 & *532 = & \pm \sI_{60} & \mbox{pentakis dodecahedron} \\ 
17~~~ & 35.2353 & 35.8078 & \times = & \pm \sC_1 &  \\
18~~~ & 34.4088 & 35.1897 & 3\times = & \pm \sC_3 &  \\
19~~~ & 33.2115 & 34.2507 & \times = & \pm \sC_1 &  \\
20~~~ & 32.7071 & 33.1584 & *222 = & \pm \sD_4 &  \\
21~~~ & 32.2161 & 32.5064 & 5\times  =& \pm \sC_5 &  \\
22~~~ & 31.8963 & 31.9834 & 2\! * \!3 = & \pm \sD_6 &  \\
23~~~ & 30.5062 & 30.9592 & \times = & \pm \sC_1 &  \\
24~~~ & 30.1628 & 30.7628 & 3\! * \!2 = & \pm \sT_{12} &  \\
25~~~ & 29.2486 & 29.7530 & 3\times = & \pm \sC_3 & \\
26~~~ & 28.7126 & 29.1948 & 2* = & \pm \sC_2 & \\
27~~~ & 28.2495 & 28.7169 & \times = & \pm \sC_1 &  \\
28~~~ & 27.8473 & 28.1480 & \times = & \pm \sC_1 & \\ \hline
\end{array}
$$
\end{table}

\section{Packing planes in $\RR^4$}
\hsp
We were naturally very interested to see what would happen when we studied
packings in $G(4,2)$, that is, packings of planes in $\RR^4$.
Here for the first time we felt we were sailing in waters where no one had been before.
(Of course, as already mentioned in Section 2, we expected that $G(4,2)$ would be special.)

Using our optimizer, we looked for the best packings we could find of $N$ planes
in $\RR^4$ (i.e. $N$ points in $G(4,2))$, for values of $N$ up to about 50.

The coordinates of the planes as found by the computer are with respect to a random coordinate frame, and so must be ``beautified'' by hand.
Initially they look like this
$$
\begin{array}{rrrr}
-0.4909573575989161 & -0.5698930951299707 & -0.6378400314190236 & 0.1653566674254520 \\
-0.7841182881275834 & 0.0542768735416466 & 0.4433546762436948 & -0.4308702383261817 \\
-0.0989775328964588 & -0.2213935599532070 & -0.9654594563867928 & -0.0952700250146647 \\
\multicolumn{1}{c}{\cdots} &
\multicolumn{1}{c}{\cdots} &
\multicolumn{1}{c}{\cdots} &
\multicolumn{1}{c}{\cdots}
\end{array}
$$
That is, they have no obvious structure.
The reason they are called ``planes'' is that they are indeed ``plain'' (at this point I held up a sheet of perfectly blank white paper;
the audience laughed).
The fun comes in trying to understand the computer output.

For $N=2$ planes, the best packings for both definitions of distance are the same:
take two orthogonal planes,
say
$$
\begin{array}{llll}
+ & 0 & 0 & 0 \\
0 & + & 0 & 0
\end{array}
\qquad {\rm and} \qquad
\begin{array}{llll}
0 & 0 & + & 0 \\
0 & 0 & 0 & +
\end{array} ~,
$$
with principal angles $\pi /2$,
so that $\pi /2$, $d_c^2 =2$, $d_g^2 = \pi^2 /2$.
(We abbreviate $+1$ and $-1$ by $+$ and $-$ respectively.)

For $N=3$ planes there are different answers for the two distances.
For the geodesic distance the best packing consists (for example)
of the planes
$$
\begin{array}{llll}
+ & - & 0 & 0 \\
0 & 0 & 0 & +
\end{array}~, \qquad
\begin{array}{llll}
+ & 0 & - & 0 \\
0 & + & 0 & 0
\end{array}~, \qquad
\begin{array}{llll}
+ & 0 & 0 & - \\
0 & 0 & + & 0
\end{array}~,
$$
and has $d_g^2 = 5 \pi^2 /18$, $d_c^2 = 1.25$.
For chordal distance the best packing consists (for example)
of the planes
$$
\begin{array}{rrrr}
1 & 0 & r & .5 \\
0 & 1 & -.5 & r
\end{array}~, \qquad
\begin{array}{rrrr}
1 & 0 & -r & .5 \\
0 & 1 & -.5 & -r
\end{array}~, \qquad
\begin{array}{rrrr}
1 & 0 & 0 & -1 \\
0 & 1 & 1 & 0
\end{array}~,
$$
and has $d_c^2 = 1.5$, $d_g^2 = 2 \pi^2 /9$.
A more geometrical description of these two arrangements will be given below.

I postpone discussion of $N=4$ and 5,
and consider $N=6$ next.
Here the same arrangement appeared to be optimal for both distances.
In this packing the principal angles between {\em any} two of the six planes
are $\pi /2$ and ${\rm arccos}~ 2/ \sqrt{5}$, so the 6 planes lie at the vertices
of a regular simplex in $G(4,2)$.
This suggested that the arrangement should somehow be related to the icosahedron, since the angles
between any two of the six diameters of the icosahedron are equal
(i.e. these six diameters form a regular simplex in $G(3,1)$).
We soon realized that there is a classical theorem which explains this,
and gives a simple way of describing $G(4,2)$.

We remind the reader that any element $\alpha$ of $SO(4)$ may be represented
as
$$\alpha : x \mapsto \bar{\ell} xr ~,$$
where $x=x_0 + x_1 i + x_2 j + x_3 k$ represents a point on $S^3$ and
$\ell$, $r$ are unit quaternions \cite{DV64}.
The pair $- \ell$, $-r$ represent the same $\alpha$.
The correspondence between $\alpha$ and $\pm ( \ell , r)$ is one-to-one.

Given a plane $P \in G(4,2)$, let $\alpha$ be the element of $SO(4)$ that fixes $P$ and negates the points of the orthogonal plane $P^\perp$.
Then $\alpha^2 =1$, and for this $\alpha$,
it is easy to see that
$\ell = \ell_1 i + \ell_2 j + \ell_3 k$ and $r= r_1 i + r_2 j + r_3 k$ are purely
imaginary unit quaternions.
This establishes the following result (which can be found for example
in Leichtweiss \cite{Lei61}).
\begin{theorem}\label{th1}
A plane $P \in G (4,2)$ is represented by a pair $(\ell,r) \in S^2 \times S^2$, with $(- \ell , -r)$ representing the same plane.
\end{theorem}

Given two planes $P,Q \in G(4,2)$, represented by $\pm ( \ell , r)$, $\pm ( \ell ' , r' )$, respectively, the principal angles $\theta_1$, $\theta_2$ between them may be found as follows.
Let $\phi$ (resp. $\psi$) be the angle between $\ell$ and $\ell '$ (resp. $r$ and $r'$), with $0 \le \phi , \psi \le \pi$.
If $\phi + \psi > \pi$, replace $\phi$ by $\pi - \phi$ and $\psi$ by
$\pi - \psi$,
so that $0 \le \phi + \psi \le \pi$, with $\phi \le \psi$ (say).
Then
$$
\theta_1 , \theta_2 = \frac{\psi \pm \phi}{2} ~,
\quad
d_g^2 (P,Q)  =  \frac{\psi^2 + \phi^2}{2} ~, \quad
d_c^2 (P,Q)  =  1- \cos \psi ~ \cos \phi ~.
$$

A set $S = \{ P_1, \dd, P_N \} \subseteq G(4,2)$ is thus represented by
a ``binocular code'' consisting of a set of pairs
$\pm ( \ell_i , r_i) \in S^2 \times S^2$.
We call the list of $2N$ points $\pm \ell_i$ (they need not be distinct)
the ``left code'' corresponding to $S$, and the points $\pm r_i$ the ``right code''.
Conversely, given two multisets $L \subseteq S^2$, $R \subseteq S^2$,
each of size $2N$ and
closed under negation, and a bijection or ``matching'' $f$ between
them that satisfies $f( - \ell ) = -f ( \ell )$, $\ell \in L$, we obtain a set of $N$ planes in $G(4,2)$.

\begin{figure}[htb]
\centerline{\psfig{file=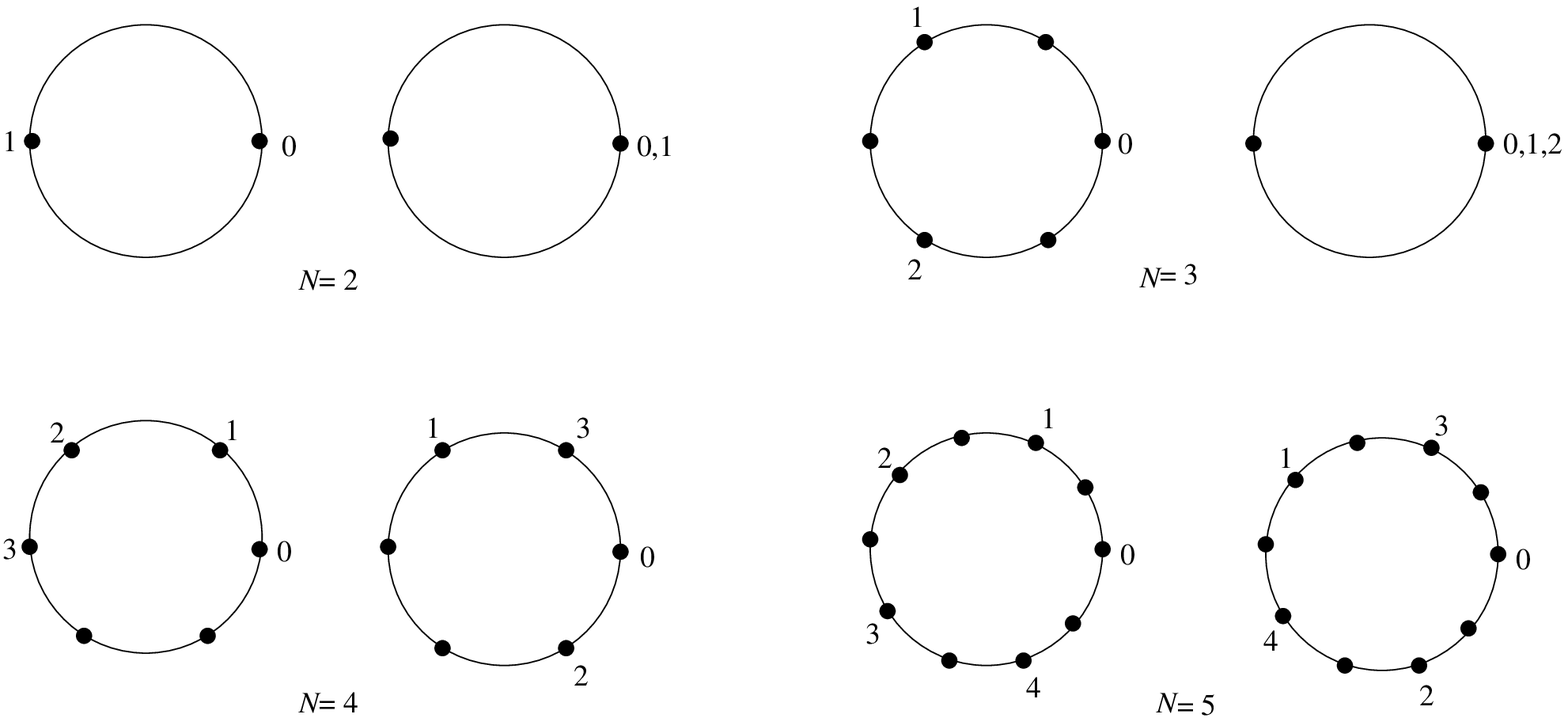,width=5.75in}}

\caption{Binocular codes describing best packings of $N=2 , \dd , 5$ planes in $G(4,2)$ for chordal distance.}
\label{fg11}
\end{figure}

\begin{figure}[htb]
\centerline{\psfig{file=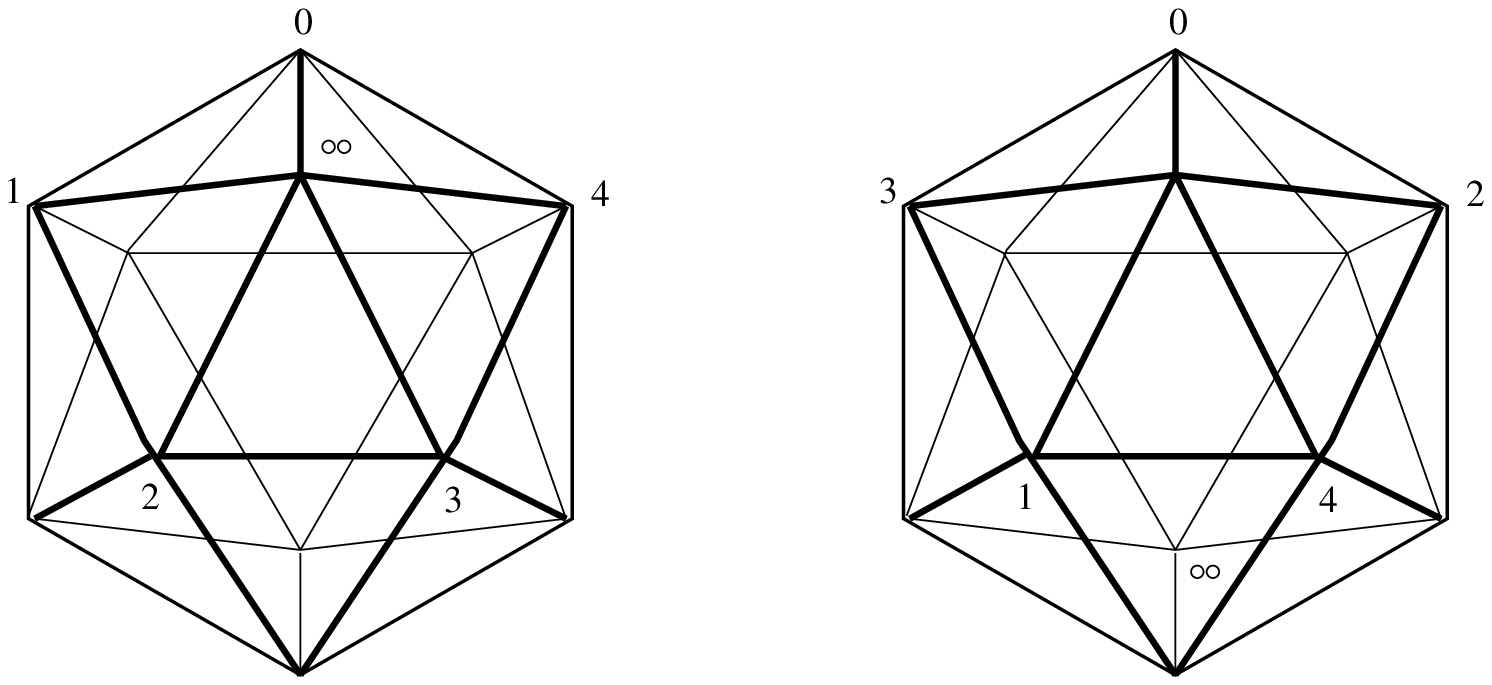,width=5in}}

\caption{Best packing of 6 planes in $G(4,2)$ with respect to both
metrics.
The left and right codes comprise the vertices of an icosahedron.
Adjacent vertices in one code are matched with non-adjacent vertices in the other code.}
\label{fg12}
\end{figure}

The binocular codes for the $d_c$-optimal packings of $N=2 , \dd , 6$ planes are shown in Figs.~\ref{fg11}, \ref{fg12}.
Except for $N=3$, the left and right codes are identical.
Matching points from the left and right codes are labeled with the same symbol.
For $N \le 5$ the points lie in the equatorial plane, and for $N \le 4$ there are repeated points.
The points lie on regular figures, except for $N=4$ where the points are
$\pm (1,0,0)$, $\pm \left( \frac{1}{\sqrt{3}} , \pm \sqrt{\frac{2}{3}} , 0 \right)$.

For $N=6$ the left and right codes consist of the 12 vertices of an icosahedron
(Fig.~\ref{fg12}).
Let these be the points
$$
\la (0, 
\pm 1, \pm \tau ) ,~
\la ( \pm \tau , 0, \pm 1 ),~
\la ( \pm 1, \pm \tau , 0) ~,
$$
where $\la = 1 / \sqrt{ \tau +2}$.
The matching is obtained by mapping each point to its algebraic conjugate
(i.e. replacing $\sqrt{5}$ by $- \sqrt{5}$), and rescaling so the points
again lie on a unit sphere.
As already mentioned, the principal angles between each pair of these
planes are $\arcsin~ 1/\sqrt{5}$ and $\pi/2$, so $d_c^2 = 6/5$, $d_g^2 = 2.6824$.

Here is another very nice packing, also found by the algorithm:
for 18 planes in $\RR^4$, use the binocular code consisting of the pairs
$(l,r)$, where $l$ and $r$ range over the vertices of a regular octahedron.
This has $d_c^2 =1$ (which is optimal, by the Corollary below) and
$d_g^2 = \pi^2 /8$ (not optimal).

Table~\ref{TD4} gives the values of $d_c^2$ for the best packings we have found of $N \le 50$ planes in $\RR^4$.

\begin{table}[htb]
\caption{Values of $d_c^2$ for best packings found of $N \le 50$ planes in $\RR^4$.}
\label{TD4}

$$
\begin{array}{ccccccccccc}
N & 3 & 4 & 5 & 6 & 7 & 8 & 9 & 10 & 11 & 12 \\
d_c^2 & 1.5000 & 1.3333 & 1.2500 & 1.2000 & 1.1667 & 1.1429 & 1.1231 & 1.1111 & 1.0000 & 1.0000  \\ [+.2in]
N & 13 & 14 & 15 & 16 & 17 & 18 & 19 & 20 & 21 & 22 \\
d_c^2 & 1.0000 & 1.0000 & 1.0000 & 1.0000 & 1.0000 & 1.0000 & 0.9091 & 0.9091 & 0.8684 & 0.8629  \\ [+.2in]
N & 23 & 24 & 25 & 26 & 27 & 28 & 29 & 30 & 31 & 32 \\
d_c^2 & 0.8451 & 0.8372 & 0.8275 & 0.8144 & 0.8056 & 0.8005 & 0.7889 & 0.7809 & 0.7760 & 0.7691  \\ [+.2in]
N & 33 & 34 & 35 & 36 & 37 & 38 & 39 & 40 & 41 & 42 \\
d_c^2 & 0.7592 & 0.7549 & 0.7489 & 0.7477 & 0.7286 & 0.7198 & 0.7095 & 0.7066 & 0.6992 & 0.6948  \\ [+.2in]
N & 43 & 44 & 45 & 46 & 47 & 48 & 49 & 50 \\
d_c^2 & 0.6844 & 0.6831 & 0.6809 & 0.6793 & 0.6732 & 0.6667 & 0.6667 & 0.6667
\end{array}
$$
\end{table}

Many of the other entries in Table~\ref{TD4} are also very beautiful:
see \cite{grass}.
On the other hand, with the exception of small values of $N$ (6, for example), the best packings with respect to the geodesic distance are much uglier.
That is, the packings are less symmetric, i.e. the orders of the automorphism groups of the best geodesic packings are almost always small,
whereas those of the best chordal-distance packings are occasionally quite large (we will see an example in the next section).

\section{Packing $n$-spaces in $\RR^n$}
\hsp
We also applied our algorithm to look for examples of packings in $\RR (m,n)$, and amassed a large number of examples.
In looking for an analogue of Theorem~1 we collected
the following pieces of evidence:

(i)~There is a packing of 10 planes in $\RR^4$ that forms a regular simplex.

(ii)~The 18-plane packing described in the previous section has the structure
of a regular orthoplex (a generalized octahedron, or cross-polytope).

(iii)~We wrote a computer program to determine the lowest
dimensions into which our library of packings
in $G(m,n)$ could be isometrically embedded.
More precisely, for a given set of $N$ points in $G(m,n)$,
we searched for the smallest dimension $D$ such that
there are $N$ points in $\RR^D$
whose Euclidean distances coincide with the chordal
distances between the points.

The results were a surprise:
it appeared that $G(m,n)$ with chordal distance could be isometrically embedded
into $\RR^D$, for $D= {\binom{m+1}{2}} -1$, independent of $n$.
Furthermore the points representing elements of $G(m,n)$ were observed
to lie on a sphere of radius $\sqrt{n(m-n)/2m}$ in $\RR^D$.

These three facts, and a number of other pieces of evidence mentioned in
\cite{grass}, led us to the main theorem of that paper.
The key idea is simply to
associate to each $P \in G(m,n)$ the orthogonal projection map from
$\RR^m$ to $P$.
If $A$ is a generator matrix for $P$ whose rows are orthogonal unit vectors, then the projection is represented by the matrix
\beql{Eq29}
\sP = A^{tr} A~.
\eeq
$\sP$ is an $m \times m$ symmetric idempotent matrix, which is independent of the particular
orthonormal generator matrix used to define it.
Changing to a different coordinate frame in $\RR^m$ has the effect of conjugating
$\sP$ by an element of $O(m)$.
With the help of \eqn{Eq3}, we see that
$$
{\rm trace} ~ \sP = n ~.
$$
Thus $\sP$ lies in a space of dimension
${\binom{m+1}{2}} -1$.

Let $\| ~ \|$ denote the $L_2$-norm of a matrix:
if $M= (M_{ij} )$, $1 \le i,j \le m$,
$$\| M \| = \sqrt{\sum_{i=1}^m \sum_{j=1}^m M_{ij}^2}
= \sqrt{{\rm trace} ~M^{tr} M} ~.
$$
For $P,Q \in G(m,n)$, with orthonormal generator matrices $A$, $B$,
and principal angles $\theta_1, \dd, \theta_n$, an elementary
calculation using \eqn{Eq3}, \eqn{Eq4} shows that
\begin{eqnarray}
\label{Eq32}
d_c^2 (P,Q) & = & n- ( \cos^2 \theta_1 + \cdots + \cos^2 \theta_n ) \nonumber \\
& = & n - {\rm trace} ~ A^{tr} A B^{tr} B \nonumber \\
& = & \frac{1}{2} \| \sP - \sQ \|^2 ~,
\end{eqnarray}
where $\sP$, $\sQ$ are the corresponding projection matrices.

Note that if we define the ``de-traced'' matrix
$\bar{\sP} = \sP - \frac{n}{m} I_m$, then ${\rm trace}~ \bar{\sP} =0$, and
$\| \bar{\sP} \|^2 = \frac{n(m-n)}{n}$.
We have thus established the following theorem.
\begin{theorem}\label{th2}
The representation of $n$-spaces $P\in G(m,n)$ by their
projection matrices $\bar{\sP}$ gives an isometric
embedding of $G(m,n)$ into a sphere of radius $\sqrt{n(m-n)/n}$ in $\RR^D$,
$D= {\binom{m+1}{2}} -1$, with $d_c (P,Q) = \frac{1}{\sqrt{2}} \| \bar{\sP} - \bar{\sQ} \|$.
\end{theorem}

Thus chordal distance between spaces is $1/ \sqrt{2}$ times the straight-line distance between the projection matrices (which explains our name for this metric).
The geodesic distance between the spaces is $1/ \sqrt{2}$ times the geodesic distance between
the projection matrices measured along the sphere in $\RR^D$.

Incidentally
the Pl\"{u}cker embedding, in which members of $G(m,n)$ are
represented by points in projective space of dimension
${\binom{m}{n}} -1$, does not give a way to realize either $d_c$ or $d_g$ as Euclidean distance.
Note also that the dimension of the Pl\"{u}cker embedding is in general much larger than the dimension of our embedding.

Since we have embedded $G(m,n)$ into a sphere of radius $\sqrt{n(m-n)/m}$ in $\RR^D$,
we can apply the Rankin bounds for spherical codes \cite{Ran55}, and deduce:
\paragraph{Corollary.}
{\em
(i)~The simplex bound:
for a packing of $N$ $n$-spaces in $\RR^m$,
\beql{Eq40}
d_c^2 \le \frac{n(m-n)}{m} \cdot \frac{N}{N-1} ~.
\eeq
Equality requires $N \le D+1 = {\binom{m+1}{2}}$, and occurs if and only if the $N$ points in $\RR^D$
corresponding to the $n$-spaces form a regular `equatorial' simplex.

(ii)~The orthoplex bound: for $N > {\binom{m+1}{2}}$,
\beql{Eq41}
d_c^2 \le \frac{n(m-n)}{m} ~.
\eeq
Equality requires $N \le 2D = (m-1) (m+2)$, and occurs if the $N$ points form a
subset of the $2D$ vertices of a regular orthoplex.
If $N=2D$ this condition is also necessary.
}

\vspace*{+.1in}
The corollary allows us to establish the optimality of hundreds of our packings.

The case of subspaces of dimension $n$ in $\RR^{2n}$ or $\RR^{2n+1}$ is especially interesting.
The largest possible arrangements of subspaces that could achieve the two bounds are:
\beql{Eq42}
\begin{array}{rrr@{~}lr@{~}l}
\multicolumn{1}{c}{m} & \multicolumn{1}{c}{n} &
\multicolumn{2}{c}{N({\rm simplex})} & \multicolumn{2}{c}{N({\rm orthoplex})} \\ [+.1in]
2 & 1 & 3 & \surd & 4 & \surd \\
3 & 1 & 6 & \surd & 10 & \\
4 & 2 & 10 & \surd & 18 & \surd \\
5 & 2 & 15  && 28 & \\
6 & 3 & 21 && 40 \\
7 & 3 & 28 & \surd & 54 & \\
8 & 4 & 36 & & 70 & \surd \\
\cdot & \cdot & \cdot & & \cdot
\end{array}
\eeq
Checks indicate that such a packing exists.
It is known that the orthoplex bound cannot be achieved by 10 lines
in $G(3,1)$,
while the other cases are undecided.
Our computer experiments strongly suggest that no set of 15 planes meets
the simplex bound in $G(5,2)$.


\paragraph{70 4-spaces in $\RR^8$.}
On the other hand, it is possible to find
packings of 70 points in $G(8,4)$ meeting the bound
\eqn{Eq41}.
With a considerable amount of effort we determined several
examples, of which the following is the most symmetrical.
Let the coordinates be labeled $\infty , 0,1, \dd, 6$, and take two 4-spaces generated by the vectors
\beql{Eq55}
\begin{array}{llll}
\{10000000, & 01000000, & 00100000, & 00001000 \} \\ [+.05in]
\{11000000, & 00101000, & 00010001, & 00000110 \}
\end{array} ~,
\eeq
respectively.
We obtain 70 4-spaces from these by negating any even number of coordinates, and/or
applying the permutations $(0123456)$,
$(\infty 0) (16) (23) (45)$ and
$(124)(365)$.
The principal angles are $0,0, \frac{\pi}{2} , \frac{\pi}{2}$;
$\frac{\pi}{4}, \frac{\pi}{4} , \frac{\pi}{4} , \frac{\pi}{4}$;
or
$\frac{\pi}{2}, \frac{\pi}{2} , \frac{\pi}{2} , \frac{\pi}{2}$, so
$d_c^2 = 2$, $d_g^2 = \pi^2 /4$ (this is not even a local optimum with respect
to geodesic distance).

In fact, the automorphism group $G$ of this packing acts transitively on the 70 subspaces, so we can obtain the packing by taking (say) the first of the
spaces in \eqn{Eq55} and letting the group act.
The group has structure $2^8 \sA_8$ and order 5160960 (where $\sA_8$ is the
alternating group of order 8).

\section{The Miraculous Enters}
\hsp
A few weeks after that packing was discovered, my colleague Peter Shor,
who was studying fault-tolerant quantum computation, asked me about the
best way to investigate a certain group of $8 \times 8$ orthogonal
matrices.
I replied by citing the computer algebra system MAGMA
\cite{Mag1}, \cite{Mag2}, \cite{Mag3},
and gave as an illustration the MAGMA commands needed to specify the group $G$
of the packing of 70 4-spaces in $\RR^8$.
To our astonishment the two groups turned out (apart from a minor re-ordering
of the coordinates) to be identical
(not just isomorphic)!

We then discovered that this group was a member of an infinite family of
groups that played a central role in a joint paper \cite{CCKS96} written by
another colleague, A. R. Calderbank.
This is a certain family of Clifford groups
(the name is due to Wall
\cite{SS4}, \cite{SS5}, \cite{Wall}), which may be constructed as follows.

The starting point is the standard method of associating a finite
orthogonal space to an extraspecial 2-group, as described for
example in \cite{Asch}, Theorem~23.10,
or \cite{Huppert}, Theorem~13.8.
The end result will be the construction of various packings of
$n$-spaces in
a parent space $V = \RR^m$, where $m = 2^i$.
As basis vectors for $V$ we use $e_u$, $u \in U = \FF_2^{^i}$.
The constructions will involve certain subgroups of the
real orthogonal group $\sO = O(V, \RR )$.

For $a,b \in U$ we define transformations
$X(a) \in \sO$, $Y (b) \in \sO$ by
$$
X(a): e_u \to e_{u+a} ~, \quad
Y(b): e_u \to (-1)^{b \cdot u} e_u ~, ~
u \in U~,
$$
where the dot indicates the usual inner product in $U$.
Then $X = \langle X(a) : a \in U \rangle$,
$Y = \langle Y(b): b \in U \rangle$ are elementary
abelian subgroups of $\sO$ of order $2^i$, and
$E = \langle X,Y \rangle \subset \sO$ is an
extraspecial 2-group$^3$\footnotetext[3]{\cite{Huppert}, p.~349} of order
$2^{2i+1}$
(\cite{CCKS96}, Lemma~2.1).
The elements of $E$ have the form $\pm X(a) Y (b)$,
$a,b \in U$, and satisfy
$$
Y(b) X(a) = (-1)^{a \cdot b}
X(a)Y(b)~,
$$
$$
(-1)^s X(a) Y(b) (-1)^{s'}
X(a') Y(b') =
(-1)^{a' \cdot b + s+s'}
X(a+a') Y(b+b') ~.
$$

The center $\Xi (E)$ of $E$
is $\{ \pm I\}$, and $\bar{E} = E/ \Xi (E)$ is an elementary
abelian group of order $2^{2i}$ whose elements
can be denoted by
$\bar{X}(a) \bar{Y}(b)$,
$a,b \in U$, where we are using the bar $\bar{~~}$ for images under the
homomorphism from $E$ to $\bar{E}$.
As in \cite{Asch}, Theorem~23.10 we define a
quadratic form $Q: \bar{E} \to \FF_2$ by
$$
Q( \bar{g} ) = \left\{
\begin{array}{ll}
0 & \mbox{if $g^2 = + I$} \\
1 & \mbox{if $g^2 = - I$}
\end{array}
\right.
$$
for $\bar{g}\in \bar{E}$, where $g \in E$ is any preimage of $\bar{g}$,
and so $Q( \bar{X} (a) \bar{Y} (b)) = a \cdot b$.

The associated alternating bilinear form $B: \bar{E} \times \bar{E} \to \FF_2$
is given by
$$
B( \bar{g}_1, \bar{g}_2 ) = Q( \bar{g}_1 +
\bar{g}_2 ) + Q( \bar{g}_1 ) + Q( \bar{g}_2 )~,
$$
for $\bar{g}_1 , \bar{g}_2 \in \bar{E}$, and so
\beql{eq5}
B( \bar{X} (a) \bar{Y} (b)~,~ \bar{X} (a')
\bar{Y} (b')) = a \cdot b' + a' \cdot b~.
\eeq
Then $(\bar{E},Q)$ is an orthogonal vector space of type $\Om^+ (2i,2 )$
and maximal Witt index (cf. \cite{Dieu}).

The Clifford group $L$ that we need is the
normalizer of $E$ in $\sO$.
This has order
$$
2^{i^2 + i +2}~~ (2^i -1)~~ \Pi_{j=1}^{i-1} (4^j -1 )
$$
(cf. \cite{CCKS96}, Section~2).
For $i=3$ the order is 5160960:
this is the group mentioned at the end of the last section.
$L$ is generated by $E$, all permutation matrices
$G(A,a) \in \sO : e_u \to e_{Au+a} $,
$u \in U$, where
$A$ is an invertible $i \times i$ matrix over $\FF_2$ and $a \in U$, and the further matrix
$H = (H_{u,v} ) $,
$H_{u,v} = 2^{- i/2} (-1)^{u \cdot v}$,
$u,v \in U$.

The group $L$ acts on $E$ by conjugation, fixing the center, and
so also acts on $\bar{E}$.
In fact $L$ acts on $\bar{E}$ as the orthogonal
group $O^+ (2i,2)$
(\cite{CCKS96}, Lemma~2.14).

This Clifford group $L$ has arisen in several different contexts,
providing a link between the the problem of packing in Grassmannian spaces,
the Barnes-Wall lattices (see \cite{SS4},
\cite{SS5}, \cite{grass2}, \cite{Wall}), the construction of
orthogonal spreads and Kerdock sets \cite{CCKS96}, and the
construction of quantum error-correcting codes \cite{BDSW},
\cite{qc2}.
It also occurs in several purely group-theoretic contexts -- see
\cite{CCKS96} for references.

The connection with quantum computing arises because if certain
conditions are satisfied the invariant
subspaces mentioned in Theorem~1 form good
quantum-error-correcting codes \cite{CRSS96}, \cite{qc2}.
\section{The construction from totally singular subspaces}
\hsp
We give one theorem, taken from \cite{grass3}, to illustrate how the Clifford group may be used to construct packings in Grassmannian spaces.

A subspace $\bar{S} \subseteq \bar{E}$ is {\em totally singular}
if $Q( \bar{g}) = 0$ for all
$\bar{g} \in \bar{S}$.
Then $\dim \bar{S} \leq i$,
and if $\dim \bar{S} = i$ then $\bar{S}$ is
{\em maximally totally singular}.
It follows from (\ref{eq5}) that the preimage
$T \subseteq E$ of a maximally totally singular space $\bar{T}$ is an
abelian subgroup of $E$, of order $2^{i+1}$.
$T$ contains $-I$, and has $2^{i+1}$ linear characters,
associated with $2^i$ mutually perpendicular
1-dimensional invariant subspaces forming a
coordinate frame $\sF (T) \subset V$
(\cite{CCKS96}, Lemma~3.3).

Since $L$ acts as $O^+ (2i,2)$ on $\bar{E}$,
$L$ takes
any ordered pair of maximally
totally singular subspaces that meet in $\{0\}$ to $X$ and $Y$ respectively.
The corresponding coordinate frames in $V$ are
\beql{eq6}
\sF (X) = \{ e_v^\ast = \df{1}{2^{i/2}} \sum_{u \in U}
(-1)^{u \cdot v} e_u : v \in U \}
\eeq
and
\beql{eq7}
\sF (Y) = \{e_u : u \in U \}~,
\eeq
respectively.

If $\bar{S} \subseteq \bar{T}$ has dimension $k$, its
preimage $S \subseteq E$ has $2^{k+1}$ linear characters, and
$2^k$ distinct invariant subspaces, each of which is spanned by
$2^{i-k}$ of the vectors in $\sF (T)$.

The following theorem produces many good Grassmannian packings.
For the proof see \cite{grass3}.
\begin{theorem}\label{th3}
Given $k$, with $0 \leq k \leq i-1$, the set of all invariant subspaces of the preimages $S$
of all $(i-k)$-dimensional totally singular subspaces $\bar{S}$ of $\bar{E}$
is a packing of $N$ planes in $G(2^i , 2^k )$ with minimal distance
$d = 2^{(k-1)/2}$, where
$$
N = 2^{i-k} \left[
\begin{array}{l}
i \\
k
\end{array}
\right] \prod_{j=k}^{i-1} (2^j +1) ~,
$$
{\em and}
$$
\left[
\begin{array}{l}
i \\
k
\end{array}
\right] = \df{(2^i -1) \ldots (2^{i-k+1} -1)}{(2^k -1) \ldots (2-1 )}
$$
is a Gaussian binomial coefficient.
\end{theorem}
\paragraph{\bf Examples.}
Taking $k=0$ in the theorem we obtain a packing of
$$
(2+2)(2^2 +2) \ldots (2^i +2)
$$
lines in $G(2^i , 1)$ with minimal angle $\pi/4$
(as in \cite{grass2}).
These are the lines defined by the minimal vectors in the
$2^i$-dimensional Barnes-Wall lattice together with their images
under $H$
(cf. \cite{SPLAG}, p.~151).

With $i = 2,~k=1$ and $i = 3,~k=2$
we obtain two important special cases:
18~points in $G(4,2)$ and 70~points in $G(8,4)$.
More generally, when $k = i -1$
we obtain the packing of
$$
f(i) = 2(2^i -1)(2^{i-1} +1 )
$$
points in $G(2^i , 2^{i-1})$ with $d^2 = 2^{i-2}$ that
is the main result of \cite{grass2}.
These packings meet the orthoplex bound of (\ref{Eq41}) and are
therefore optimal.
An explicit recursive construction for the special case
$k = i-1$ is given in \cite{grass2}.

For $k =1$ and $k = i-2$ we obtain two further sequences of
packings whose existence was conjectured in \cite{grass2}.

The construction given in the theorem can be restated in an equivalent but more
explicit way as follows.
Let $P_0$ be the $2^k$-dimensional space spanned by the
coordinate vectors $e_u$, where $u \in U$ is of the form
$00 \ldots 0 \ast \ldots \ast$,
with $i-k$ initial zeros.
Then the packing consists of all the images of $P_0$ under the group $L$.

The paper \cite{grass3} also contains many other Grassmannian packings.

Space does not permit any discussion of the application of the Clifford
group to construct quantum error correcting codes:
for this see \cite{qc2}.

We conclude by mentioning the present status of the existence of the packings
listed in Eq. (\ref{Eq42}).
No packing of
10 points in $G(3,1)$ meeting
the orthoplex bound
can exist, and we conjecture no packing of 15 points in $G(5,2)$ meeting the simplex bound can exist.
We have found an infinite family of packings of $p(p+1)/2$ points in $G(p,(p-1)/2)$ meeting the simplex bound, generalizing the packing of 28 points
in $G(7,3)$ \cite{grass3}.
As mentioned above, we have also found an infinite family of packings of $m^2 + m-2$ points in $G(2^m, 2^{m-1} )$ meeting the orthoplex bound,
generalizing the packings of 18 points in $G(4,2)$ and 80 points
in $G(8,4)$.
We would very much like to know if packings of 21 or 40 points in $G(6,3)$,
54 points in $G(7,3)$ or 36 points in $G(8,4)$ meeting the bounds exist.
\paragraph{Acknowledgements.}
I should like to express my thanks to my coauthors (Rob Calderbank, John Conway, Ron Hardin, Eric Rains and Peter Shor) of the papers to which this
article serves as introduction.
\clearpage

\end{document}

%% file: tent_text.tex
\begin{picture}(0,0)%
\epsfig{file=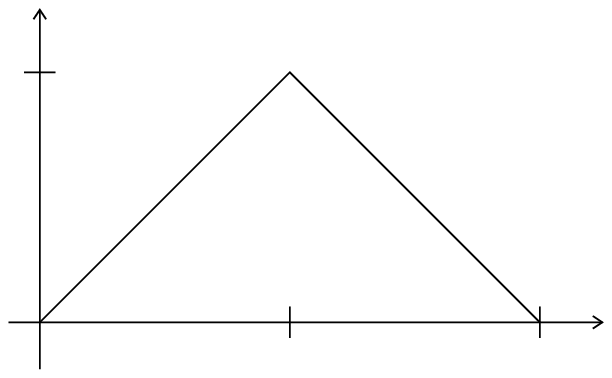}%
\end{picture}%
\setlength{\unitlength}{0.01250000in}%
\begingroup\makeatletter\ifx\SetFigFont\undefined
\def\x#1#2#3#4#5#6#7\relax{\def\x{#1#2#3#4#5#6}}%
\expandafter\x\fmtname xxxxxx\relax \def\y{splain}%
\ifx\x\y   
\gdef\SetFigFont#1#2#3{%
  \ifnum #1<17\tiny\else \ifnum #1<20\small\else
  \ifnum #1<24\normalsize\else \ifnum #1<29\large\else
  \ifnum #1<34\Large\else \ifnum #1<41\LARGE\else
     \huge\fi\fi\fi\fi\fi\fi
  \csname #3\endcsname}%
\else
\gdef\SetFigFont#1#2#3{\begingroup
  \count@#1\relax \ifnum 25<\count@\count@25\fi
  \def\x{\endgroup\@setsize\SetFigFont{#2pt}}%
  \expandafter\x
    \csname \romannumeral\the\count@ pt\expandafter\endcsname
    \csname @\romannumeral\the\count@ pt\endcsname
  \csname #3\endcsname}%
\fi
\fi\endgroup
\begin{picture}(244,120)(58,542)
\put(200,545){\makebox(0,0)[b]{\smash{\SetFigFont{12}{14.4}{rm}$\pi /2$}}}
\put(280,544){\makebox(0,0)[b]{\smash{\SetFigFont{12}{14.4}{rm}$\pi$}}}
\put(240,545){\makebox(0,0)[b]{\smash{\SetFigFont{12}{14.4}{rm}$\phi$}}}
\put(113,596){\makebox(0,0)[rb]{\smash{\SetFigFont{12}{14.4}{rm}$\theta_1$}}}
\put(127,635){\makebox(0,0)[lb]{\smash{\SetFigFont{12}{14.4}{rm}$\pi/2$}}}
\put(125,547){\makebox(0,0)[lb]{\smash{\SetFigFont{12}{14.4}{rm}$0$}}}
\end{picture}

%% file: yama.bbl
\begin{thebibliography}{99}
\bibitem{Asch}
M. Aschbacher,
{\em Finite Group Theory},
Cambridge Univ. Press, 1986.

\bibitem{Asi85}
D. Asimov,
The Grand Tour --- a tool for viewing multidimensional data,
{\em SIAM J. Sci. Stat. Comput.}, {\bf 6} (1985), 128--143.

\bibitem{BDSW}
C. H. Bennett, D. DiVincenzo, J. A. Smolin and W. K. Wootters,
``Mixed state entanglement and quantum error correction,''
{\em Phys. Rev. A}, {\bf 54} (1996), 3824--3851;
also LANL e-print quant-ph/9604024.

\bibitem{SS4}
B. Bolt, T. G. Room and G. E. Wall,
On Clifford collineation, transform and similarity groups I,
{\em J. Australian Math. Soc.}, {\bf 2} (1961), 60--79.
 
\bibitem{SS5}
B. Bolt, T. G. Room and G. E. Wall,
On Clifford collineation, transform and similarity groups II,
{\em J. Australian Math. Soc.}, {\bf 2} (1961), 80--96.

\bibitem{Mag1}
W. Bosma and J. Cannon,
{\em Handbook of Magma Functions},
Sydney, May~22, 1995.

\bibitem{Mag2}
W. Bosma,
J. J Cannon and G. Mathews,
Programming with algebraic structures:
Design of the Magma language, In:
M. Giesbrecht (ed.),
{\em Proceedings of the 1994 International Symposium on Symbolic and Algebraic Computation},
Association for Computing Machinery, 1994, 52--57.

\bibitem{Mag3}
W. Bosma, J. Cannon and C. Playoust,
The Magma algebra system I:
The user language,
{\em J. Symb. Comp.},
{\bf 24} (1997), 235--265.

\bibitem{CCKS96}
A. R. Calderbank, P. J. Cameron, W. M. Kantor and J. J. Seidel,
``$\ZZ_4$ Kerdock codes, orthogonal spreads, and extremal Euclidean line-sets,''
{\em Proc. London Math. Soc.}, {\bf 75} (1997), 436--480.

\bibitem{grass3}
A. R. Calderbank, R. H. Hardin, E. M. Rains, P. W. Shor and N. J. A. Sloane,
``A group-theoretic framework for the construction of packings in Grassmannian spaces,''
{\em J. Algebraic Combinatorics}, 1997 (submitted).

\bibitem{CRSS96}
A. R. Calderbank, E. M. Rains, P. W. Shor and N. J. A. Sloane,
``Quantum error correction and orthogonal geometry,''
{\em Phys. Rev. Lett.},
{\bf 78} (1997), 405--409;
also LANL e-print quant-ph/9605005.

\bibitem{qc2}
A. R. Calderbank, E. M. Rains, P. W. Shor and
N. J. A. Sloane, Quantum error correction via codes over
$GF(4)$,
{\em IEEE Trans. Inform. Theory},
{\bf 44} (1998), in press.

\bibitem{grass}
J. H. Conway, R. H. Hardin, and N. J. A. Sloane,
``Packing lines, planes, etc.:
packings in Grassmannian space,''
{\em Experimental Math.}, {\bf 5} (1996), 139--159.
See also {\bf 6} (1997), p.~175.

\bibitem{SPLAG} J.~H.~Conway and N.~J.~A.~Sloane,
{\it Sphere Packings, Lattices and Groups\/}, 3rd ed.,
Grundlehren der math. Wissenschaften 290, New York: Springer, 1998.

\bibitem{CS96}
J. H. Conway and N. J. A. Sloane,
The four-dimensional point groups --- a revised enumeration, in preparation.

\bibitem{Dieu}
J. A. Dieudonn\'e,
{\em La G\'eom\'etrie des Groupes Classiques},
Springer-Verlag, Berlin, 1971.

\bibitem{DV64}
P. Du~Val,
{\em Homographies, Quaternions and Rotations},
Oxford Univ. Press, 1964.

\bibitem{FT65}
L. Fejes T\'{o}th,
Distribution of points in the elliptic plane,
{\em Acta Math. Acad. Sci. Hungar.},
{\bf 16} (1965), 437--440.

\bibitem{GVL89}
G. H. Golub and C. F. Van~Loan,
{\em Matrix Computations},
Johns Hopkins Univ. Press, 2nd ed., 1989.

\bibitem{HS92}
R. H. Hardin and N. J. A. Sloane,
{\em Operating Manual for Gosset: A General-Purpose Program for Constructing Experimental Designs (Second Edition)}.
Statistics Research Report No.~106, Bell Laboratories,
Murray Hill, NJ, October, 1992.
Also DIMACS Technical Report 93--51,
DIMACS Center, Rutgers University, New Brunswick, NJ, August 1993.

\bibitem{HS93}
R. H. Hardin and N. J. A. Sloane,
A New Approach to the Construction of Optimal Designs,
{\em J. Statistical Planning and Inference},
{\bf 37} (1993), 339--369.

\bibitem{HSS94}
R. H. Hardin, N. J. A. Sloane and W. D. Smith,
{\em Tables of Spherical Codes in Dimensions 3, 4, 5},
Netlib archive,
{\tt http://netlib.att.com/math/sloane/packings},
1993--1995.
This table is based on our own calculations, with additional
packings provided by D. A. Kottwitz and J. Buddenhagen, and is believed
to contain the best packings known with $N \le 130$ points in dimensions
3, 4 and 5.

\bibitem{HSS96}
R. H. Hardin, N. J. A. Sloane,
and W. D. Smith,
{\em Spherical Codes}, book in preparation.

\bibitem{Me174}
R. H. Hardin and N. J. A. Sloane,
New Spherical 4-Designs,
{\em Discrete Mathematics},
{\bf 106/107} (1992), 255--264.
(Topics in Discrete Mathematics, vol. 7, "A Collection
of Contributions in Honour of Jack Van Lint", ed. P. J. Cameron
and H. C. A. van Tilborg, North-Holland, 1992.)

\bibitem{Me189}
R. H. Hardin and N. J. A. Sloane,
Expressing ($a^2 + b^2 + c^2 + d^2 )^2$
as a Sum of 23 Sixth Powers,
{\em Journal of Combinatorial Theory, Series A},
{\bf 68} (1994), 481--485.

\bibitem{Me203}
R. H. Hardin and N. J. A. Sloane,
Codes (Spherical) and Designs (Experimental),
{\em Different Aspects of Coding Theory, ed. A. R. Calderbank,
AMS Series Proceedings Symposia Applied Math.},
{\bf 50} (1995), 179--206.

\bibitem{Me204}
R. H. Hardin and N. J. A. Sloane,
McLaren's Improved Snub Cube and Other New Spherical Designs in Three Dimensions,
{\em Discrete and Computational Geometry},
{\bf 15} (1996), 429--441.

\bibitem{Huppert}
B. Huppert,
{\em Endliche Gruppen},
Springer-Verlag, Berlin, 1967.

\bibitem{Lei61}
K. Leichtweiss, Zur Riemannschen Geometrie in Grassmannschen
Mannigfaltigkeiten,
{\em Math. Zeit.},
{\bf 76} (1961), 334--366.

\bibitem{MS77}
F. J. MacWilliams and N. J. A. Sloane,
{\em The Theory of Error-Correcting Codes},
North-Holland, Amsterdam, 1977.

\bibitem{Ran55}
R. A. Rankin,
The closest packing of spherical caps in $n$ dimensions,
{\em Proc. Glasgow Math. Assoc.},
{\bf 2} (1955), 139--144.

\bibitem{Ros93}
M. Rosenfeld,
Problem No.~10293,
{\em Amer. Math. Monthly}, {\bf 100} (1993), 291.

\bibitem{Ros94}
M. Rosenfeld,
How wide can you spread your chopsticks?
{\em Congressus Numerantium}, {\bf 102} (1994), 29--31.

\bibitem{Ros96}
M. Rosenfeld,
In praise of the Gram matrix, preprint.

\bibitem{Rut45}
H. Rutishauser,
\"{U}ber Punktverteilungen auf der Kugelfl\"{a}che,
{\em Comm. Math. Helv.}, {\bf 17} (1945), 327--331.

\bibitem{SchW51}
K. Sch\"{u}tte and B. L. van~der Waerden,
Auf welcher Kugel haben 5, 6, 7, 8 order 9 Punkten min Mindestabstand Eins
Platz?,
{\em Math. Ann.}, {\bf 123} (1951), 96--124.

\bibitem{grass2}
P. W. Shor and N. J. A. Sloane,
``A family of optimal packings in Grassmannian manifolds,''
{\em J. Algebraic Combinatorics}, 1998 (to appear).

\bibitem{grass0}
N. J. A. Sloane,
{\em Grassmannian packings},
{\tt http://www.research.att.com/$\sim$njas/grass}~.

\bibitem{Stroh63}
J. Strohmajer, Uber die Verteilung von Punkten auf der Kugel,
{\em Ann. Univ. Sci. Budapest Sect. Math.},
{\bf 6} (1963), 49--53.

\bibitem{TaG83}
T. Tarnai and Z. G\'{a}sp\'{a}r,
Improved packing of equal spheres on a sphere and rigidity of its graph,
{\em Math. Proc. Camb. Phil. Soc.}, {\bf 93} (1983), 191--218.
\bibitem{Wall}
G. E. Wall,
On Clifford collineation, transform and similarity groups IV,
{\em Nagoya Math. J.}, {\bf 21} (1962), 199--222.

\bibitem{Won67}
Y.-C. Wong, Differential geometry of Grassmann manifolds,
{\em Proc. Nat. Acad. Sci. USA}, {\bf 47} (1967), 589--594.

\end{thebibliography}
